\documentclass[12pt]{amsart}
\usepackage{amssymb, amsmath, amsthm, a4wide}

\def\F{\mathbb{F}}
\def\N{\mathbb{N}}
\def\R{\mathbb{R}}

\theoremstyle{plain}

\newtheorem{thm}{Theorem}

\newtheorem{lem}[thm]{Lemma}

\newtheorem*{CombinatorialLemma}{Brouwer's fan theorem}
\theoremstyle{definition}

\theoremstyle{remark}

\title{A novel proof of the Heine-Borel theorem}

\begin{document}
\thispagestyle{empty}

\author{Matthew Macauley\qquad\qquad Brian Rabern\qquad\qquad Landon
    Rabern}

\begin{abstract}
  Every beginning real analysis student learns the classic Heine-Borel
  theorem, that the interval $[0,1]$ is compact. In this article, we
  present a proof of this result that doesn't involve the standard
  techniques such as constructing a sequence and appealing to the
  completeness of the reals. We put a metric on the space of
  infinite binary sequences and prove that compactness of this space
  follows from a simple combinatorial lemma. The Heine-Borel theorem
  is an immediate corollary.
\end{abstract}

\keywords{Binary sequences, ultrametric, model theory, Brouwer's fan theorem,
  G\"odel compactness, K\"onig's infinity lemma} 

\subjclass[2000]{54E45, 03C99, 03F03}

\maketitle

\section{The Heine-Borel theorem}

Think back to your first real analysis class. In the beginning, most
of the definitions were fairly straightforward. Open and closed sets
made sense, because of the common usage of open and closed intervals
in previous math classes. It was a bit odd that open sets could also
be closed, or that sets could be neither open nor closed. But this was
``higher math,'' so you could let that one slide. Then came the
definition of ``compact.'' Completely out of nowhere. Why anyone would
ever find themselves with an open cover, let alone try to extract a
finite subcover, was beyond your wildest dreams. As you sat there in
class trying to figure out what this really meant, the professor wrote
the following two sentences on the board, with ``Heine-Borel''
preceding one of them.\footnote{Textbooks vary as to which of these
statements is called the Heine-Borel theorem and which one is a lemma
or corollary. We will refer to the compactness of $[0,1]$ as the
Heine-Borel theorem. See McCleary (2006).}
\begin{itemize}
  \item The interval $[0,1]$ is compact.
  \item A subset of $\R^n$ is compact iff it is closed and bounded. 
\end{itemize}
You might remember what came next. From an arbitrary infinite sequence
contained in $[0,1]$, a divide-and-conquer technique to construct a
particular sequence of nested intervals, from this a subsequence of
real numbers, and then a summon to the completeness of the reals, one
of those blatantly obvious analysis facts that you had no idea how to
prove (and likely still don't know). It felt a little unsatisfying, and
almost seemed like cheating. About this time, it dawned on you that
your roommate was right: mathematicians make a living saying the
simplest things in the most difficult round-about way.

By now, you understand in ways you never could have imagined back
then, how wise your old roommate was. But you also remember what
attracted you to mathematics in the first place, those mysterious
qualities that, like the silver bell in the \emph{Polar Express},
could only be heard by a select few. Your friends shook their heads
and exchanged private smiles when you marveled at the sheer beauty of
mathematics, such as the surprising connections between seemingly
unrelated topics, and the way that a basic result could be proven in
vastly different ways. In fact, it is likely fore these reasons why
you are reading this paper right now, and it is precisely for these
reasons that drove the authors to write it. So jump aboard, and enjoy
a quick but enlightening tour of diverse topics such as ultrametrics
and model theory, and we'll drop you off back in your first real
analysis course at the classic Heine-Borel theorem.


\section{A Combinatorial Lemma}

We begin by stating a combinatorial lemma due to
Brouwer~\cite{Brouwer:27}, and two proofs of it. These proofs are
quite different; one follows from K\"onig's infinity lemma, and the
other from G\"odel compactness. We will briefly explain these concepts
for those unfamiliar with them, but note the remarkable coincidence that 
both have an umlauted `o' in their name.
\begin{CombinatorialLemma}
  \label{lemma:bitstrings}
  Let $\mathcal{B}$ be a collection of finite bitstrings (binary
  sequences) so that every infinite bitstring has an initial segment
  in $\mathcal{B}$. Then there is a finite subset
  $\mathcal{A}\subseteq\mathcal{B}$ so that every infinite bitstring
  has an initial segment in $\mathcal{A}$.
\end{CombinatorialLemma}
We pause to recall K\"onig's infinity lemma, which says that an
infinite tree where every vertex has finite degree contains an
infinite path~\cite{Konig:36}. If the infinite tree is uncountable, as
ours will be, the axiom of choice is required.
\begin{proof}[(The K\"onigian Proof)]
  Assume (to reach a contradiction) that the theorem is false.
  Recursively construct a tree $T_{\F_2}$ with the empty bitstring at
  the root so that the children of $b$ are the bitstrings $b0$ and
  $b1$. Now remove all bitstrings from $T_{\F_2}$ that have an initial
  segment in $\mathcal{B}$ to get the tree $T$.
 
  For every $n\geq 1$ there exists a length-$n$ bitstring with no
  initial segment in $\mathcal{B}$ (if every bitstring of length $n$
  had an initial segment in $\mathcal{B}$, then the bitstrings from
  $\mathcal{B}$ of length at most $n$ would work for $\mathcal{A}$).
  Thus, $T$ is infinite.

  Since every vertex of $T$ has finite degree, we may apply K\"onig's
  infinity lemma to get an infinite path through $T$ starting at the
  root.  Hence we have a sequence $y_1,y_2,y_3,\dots$, where $y_i$ is
  a length-$i$ bitstring, $y_i$ is an initial segment of $y_{i+1}$,
  and none of the $y_i$'s have initial segments in $\mathcal{B}$.  Let
  $y$ be the infinite bitstring with length-$i$ initial segment $y_i$
  for each $i$. If $y$ had an initial segment of length $n$ in
  $\mathcal{B}$, then $y_n$ would have an initial segment in
  $\mathcal{B}$, which is forbidden by construction.  Hence $y$ has no
  initial segment in $\mathcal{B}$, and this contradiction completes
  the proof.
\end{proof}

We can get a more transparent proof of the lemma using G\"odel
compactness. This is a classic result from model
theory~\cite{Kreisler:71} which says that a set of boolean clauses is
satisfiable if and only if every finite subset of clauses is
satisfiable~\cite{Goedel:30}.
\begin{proof}[(The G\"odelian Proof)]
Consider a set $\{a_1,a_2,a_3,\ldots\}$ of distinct boolean variables,
and let $K$ be the following set of clauses:
\begin{align*}
  &\{N(b) \mid b \in \mathcal{B}\}\;,
\end{align*}
where $N(b) = \neg\left[(b_1 = a_1) \wedge (b_2 = a_2) \wedge \cdots
\wedge (b_k =a_k)\right]$ for any $b = b_1b_2\cdots b_k \in
\mathcal{B}$.  Note that $N(b)$ is satisfiable if and only if $b$ is
not an initial segment of $a_1a_2a_3\cdots$.

Now, assume (to reach a contradiction) that the theorem is false.
Then for any finite $\mathcal{A}\subseteq\mathcal{B}$ there exists a
bitstring $a_1a_2a_3\cdots$ that has no initial segment in
$\mathcal{A}$.  Hence, every finite subset of clauses of $K$ is
satisfiable, and by G\"odel compactness, $K$ is satisfiable.  But by
construction, this yields an infinite bitstring $a_1a_2a_3\cdots$ with
no initial segment in $\mathcal{B}$. This contradiction completes the
proof.
\end{proof}


\section{The Bit-Metric}

Equipped with our combinatorial lemma, we resume our tour in the land
of bitstrings. Let $\F_2=\{0,1\}$, and let $\F_2^\N$ denote the set of
infinite bitstrings. Again, we write a bitstring as $a=a_1a_2\cdots$,
and call the individual $a_i$'s \emph{bits}. Define the function
$\iota$ that sends an element of $\F_2^\N$ to the corresponding number
in $[0,1]$ written in binary, by
\[
\iota\colon\F_2^\N\longrightarrow[0,1]\;,
\qquad\iota(a_1a_2a_3\cdots)=\sum_{i=1}^\infty a_i
2^{-i}=0.a_1a_2a_3\dots\;.
\]
At this point, it is easy to overlook the fact that decimals have a
few subtle but pesky properties, such as the fact that
\begin{equation*}
  \label{eq:iota}
  \iota(a_1a_2\cdots a_k1000\cdots)
  =\iota(a_1a_2\cdots a_k0111\cdots)\;.
\end{equation*}
Fortunately, $\iota$ is injective on bitstrings not of this form. With
this in mind, we say that a binary decimal representation of
$x\in[0,1)$ is in \emph{standard form} if there are infinitely many
$0$s. Since nobody would ever write infinitely many $1$s instead of
just one, when we speak of a number $x\in[0,1]$, we shall assume that
it is written in standard binary form. With this assumption, we can
define the preimage of $x\in[0,1)$ under $\iota$ to be the bitstring
with infinitely many $0$s, which is denoted by $\iota^{-1}(x)$. At
this time it should be noted that the authors aren't analysts, and
thus are prone to omit crucial but obvious details, such as what the
standard binary form of $1$ is, and how to define $\iota^{-1}(1)$.

Once you resolve these tiny missing details, we may proceed together,
and put a metric on $\F_2^\N$ by saying that two distinct bitstrings
$a$ and $b$ are a distance $\beta(a,b)=2^{-k}$ apart, where $k$ is the
last bit at which $a$ and $b$ agree. It is straightforward to show
that $(\F_2^\N,\beta)$ is an \emph{ultrametric}, and we call it the
\emph{bit-metric} on $\F_2^\N$. An ultrametric is any metric that
satisfies the strong triangle inequality:
$\beta(a,c)\leq\max\{\beta(a,b),\beta(b,c)\}$, and this gives it some
extra special properties such as:
\begin{itemize}
  \item\emph{Russian doll property of balls}: If $B_r(a)\cap
  B_r(b)\neq\emptyset$, then either $B_r(a)\subseteq B_r(b)$ or
  $B_r(a)\supseteq B_r(b)$.
\item\emph{Center of the universe property}: If $|a-b|<r$, then
  $B_r(a)=B_r(b)$.
\end{itemize}
These properties are very useful when studying $(\F_2^\N,\beta)$, and
we utilize them in papers that are much more difficult to read than
this one. However, we will not need them for Heine-Borel, but we
mention them for completeness (of the paper, not the reals). 

At this point, you might be suspecting that the map $\iota$, being so
simple, is continuous. This is indeed correct since, by definition,
$|\iota(a) - \iota(b)| \leq \beta(a, b)$ for any $a, b \in \F_2^\N$.
\begin{lem}
  \label{lemma:iota}
  The map $\iota$ is continuous under the bit-metric.
\end{lem}

You might also be suspecting that under the bit-metric, $\F_2^\N$,
being a collection of infinite sequences, is not compact. This is
incorrect.
\begin{lem}
  \label{lemma:compactness}
  $(\F_2^\N,\beta)$ is compact.
\end{lem}
\begin{proof}
  Consider an open cover $\cup_{i\in I} B_{\epsilon_i}(a_i)=\F_2^\N$
  of balls, where each $a_i\in\F_2^\N$. Let $S_i$ be the first
  $\left\lfloor\log_2(\epsilon_i^{-1})\right\rfloor+1$ bits of $a_i$ for
  $i \in I$.  Then $b\in\F_2^\N$ is in $B_{\epsilon_i}(a_i)$ if and
  only if $S_i$ is an initial segment of $b$.  Put $\mathcal{B} =
  \{S_i \mid i \in I\}$ and apply Brouwer's fan theorem to get a finite
  $\mathcal{A}\subseteq\mathcal{B}$.  By construction, the set
  \[
  \bigcup_{a_i\in\mathcal{A}} B_{\epsilon_i}(a_i)=\F_2^\N\;,
  \]
  and thus we have a finite subcover of $\F_2^\N$.
\end{proof}
Equipped with Lemmas~\ref{lemma:iota} and~\ref{lemma:compactness},
we can now present The Shortest Proof of Heine-Borel Ever.
\begin{thm}[Heine-Borel]
  The interval $[0,1]$ is compact. 
\end{thm}
\begin{proof}
  $\iota(\F_2^\N)=[0,1]$ is the continuous image of a compact set.
\end{proof}

\vspace{0.15in}

This concludes our tour, now that we have arrived back at your first
real analysis class, on that special day when you first saw the
Heine-Borel theorem proven. For those of you out there that have yet
to take real analysis, but are advanced and motivated enough to be
reading this article, pay attention. When you find yourself in an
analysis class, and the professor draws that little box at the end of
the proof of Heine-Borel, raise your hand, and inquire: \\

\noindent ``\emph{Doesn't that just follow from K\"onig's infinity lemma,
and the standard ultrametric on the space of binary sequences?}'' \\


\begin{thebibliography}{99}

\bibitem{Brouwer:27} L.~E.~J. Brouwer, \"Uber
  {D}efinitionsbereiche von {F}unktionen, \emph{Math. Annalen} \textbf{97}
  (1927), 60--75. 
  
\bibitem{Diestel:05} Reinhard Diestel, \emph{Graph Theory}. Graduate
  Texts in Mathematics, \textbf{173}, Springer-Verlag Heidelberg, New
  York, 2005.

\bibitem{Goedel:30} Kurt G\"odel, Die Vollstandigbert der Axiome
  des logischen Vunktionen Kalkuls, \emph{Monatshefte f\"ur
  Mathematik und Physik}, \textbf{87} (1930), 349--360.
  
\bibitem{Konig:36} D\'enes K\H{o}nig, \emph{Theorie der Endlichen und
  Unendlichen Graphen: Kombinatorische Topologie der
  Streckenkomplexe}, Leipzig: Akad. Verlag., 1936.
  
\bibitem{Kreisler:71} H.~Jerome Kreisler, \emph{Model Theory for Infinitary
 Logic}. North-Holland Publishing Company, Amsterdam London, 1971.

\bibitem{McCleary:06} John McCleary, \emph{A First Course in Topology:
  Continuity and Dimension}, American Mathematical Society, Student
  Mathematical Library \textbf{31}, 2006.
  
\bibitem{Semmes:07} Stephen Semmes, \emph{An Introduction to the
  Geometry of Ultrametric Spaces}, (2007). arXiv:0711.0709.
  
\bibitem{VanAllsburg:85} Chris Van Allsburg. \emph{The Polar Express}, 
  Houghton Mifflin, 1985. 

\end{thebibliography}
\end{document}